\documentclass[12 pt]{article}

\usepackage{algpseudocode}
\usepackage{color}
\usepackage{graphicx}
\usepackage{bm}
\usepackage[margin=2.45cm]{geometry}
\usepackage{amsmath}
\usepackage{amssymb}
\usepackage[numbers,super,sort&compress]{natbib}
\usepackage{setspace}
\usepackage{times}

\usepackage{abstract}

\font\myfont=cmr12 at 12pt

\usepackage{authblk}

\begin{document}

\title{\vspace{-5ex}Topology-dependent density optima for efficient simultaneous network exploration}

\author[1]{Daniel B. Wilson}
\author[1]{Ruth E. Baker}
\author[2,*]{Francis G. Woodhouse}
\affil[1]{ {\myfont Wolfson Centre for Mathematical Biology,
    Mathematical Institute, University of Oxford, Radcliffe
    Observatory Quarter, Oxford OX2 6GG, United Kingdom. }}
\affil[2]{ {\myfont Department of Applied Mathematics and Theoretical Physics, Centre for Mathematical Sciences, University of Cambridge, Wilberforce Road, Cambridge CB3 0WA, United Kingdom.} }
\affil[*]{ \myfont Corresponding author: f.g.woodhouse@damtp.cam.ac.uk}
\date{}

\maketitle

\doublespacing

\begin{abstract}
{\myfont \bf A random search process in a networked environment is governed by the time it takes to visit every node, termed the cover time. Often, a networked process does not proceed in isolation but competes with many instances of itself within the same environment. A key unanswered question is how to optimise this process: how many concurrent searchers can a topology support before the benefits of parallelism are outweighed by competition for space?
Here, we introduce the searcher-averaged parallel cover time (APCT) to quantify these economies of scale.
We show that the APCT of the networked symmetric exclusion process is optimised at a searcher density that is well predicted by the spectral gap.
Furthermore, we find that non-equilibrium processes, realised through the addition of bias, can support significantly increased density optima.
Our results suggest novel hybrid strategies of serial and parallel search for efficient information gathering in social interaction and biological transport networks.}
\end{abstract}

\clearpage


From animals foraging to T-cells hunting pathogens to proteins examining DNA, nature relies on carefully optimised random searches at many scales~\cite{Benichou_ISS_review,Optimal_RunandTumble,Coppey_ProteinDNA,Bartumeus_2009,pathogens,Levy_foraging,OptimalRWs_Nature,RW_in_Bio_JRSI,2009Li_NatPhys,2016Schwarz_PRL}.
This concept is not limited to biology: robotic self-assembly~\cite{2014Rubenstein_Science}, traffic flow management~\cite{Herty_TrafficFlow,SIAMRev_Bellomo,SIAMFugen} and computer resource allocation~\cite{2001Adamic_PRE} hinge on optimising decentralised exploration.
In these distributed processes, the searcher, be it a protein, an animal or a network token, must often visit not just one site but many locations connected in a network~\cite{Chupeau_NatPhys_CTs}. The critical measure of efficiency is then the time to visit every node of the network, called the cover time~\cite{Kahn1989_CTs}. Single-searcher cover times are known on simple networks such as linear chains, rings and other regular lattices~\cite{Kundu_ExactDistinctCommon,Majumdar_ExactCTRing,Yokoi_1990_ExactCTs,CTs_2Dlattice,CTs_Lattice}, and significant progress has been made on establishing transport statistics for more general networks~\cite{MaierPRE2017,Tejedor_PRE09,Chupeau_NatPhys_CTs,Meyer_PRE12,Palacios_03,JonassonCTs,2017Weng_PRE}, providing a means to design topologies that can be searched efficiently.
However, when multiple parallel searchers compete for space or resources, as often occurs in the examples above, how to design  optimal search strategies remains a significant open question. Our central result is that the optimal density of searchers depends heavily, yet predictably, on network topology, implying that search strategies can be made efficient by careful optimisation of topology and searcher quantity.


The exclusion process represents the most fundamental model of competition for space~\cite{LiggetBook,Ghosts_Folks,PlankSimpson_JRSI}. It has been key to understanding such diverse phenomena as cell migration~\cite{2005Chowdhury_PhysLifeRev}, molecular traffic~\cite{2017Graf_PRL,2011Neri_PRL}, surface roughening~\cite{1987Plischke_PRB} and queueing~\cite{1986Kipnis_AnnProb}. Capitalising on this breadth, we use the exclusion process to model parallel searching of a network and introduce the \emph{searcher-averaged parallel cover time} (APCT), the average per-searcher time for all searchers to visit all nodes within a network.
Strategy optimisation then demands an understanding of how both network topology and searcher density impact the APCT. Consider, for example, a scenario with as many searchers as nodes. Placing all searchers on the network simultaneously results in an infinite APCT, while a simple `serial' strategy where a single searcher is placed on the network, removed once it has visited every node and replaced with a new searcher, one at a time, is almost always inefficient. It is therefore critical to determine the optimal, or most efficient, density of parallel searchers minimising the APCT (Fig.~1(a)). Through analytic and numerical results, we find that this optimal density is heavily dependent on network topology. We demonstrate that the spectral gap, which quantifies the convergence rate of a single-searcher random walk, is a strong predictor of a network's density optimum and outperforms simpler degree-based network statistics, as measured by mutual information. We provide strategies for optimal deployment of hybrid series--parallel searches allowing for construction of efficiently-explored networks. We broaden to non-equilibrium processes by generalising to flux-conserving asymmetric exclusion processes, finding a remarkable non-monotonic relationship between density optima and the spectral gap. Our work provides an accessible route into the design of optimal search strategies in complex environments involving equilibrium and non-equilibrium processes.


\section*{Results}
\subsection*{Average parallel cover time and optimal density}

\noindent
We employ a symmetric exclusion process of $M$ parallel searchers on a network with $N$ vertices (where ${M<N}$). First, an initial configuration of searchers is generated uniformly at random with each searcher occupying its own node. The searchers then perform mutually-excluding continuous-time random walks (CTRWs) with i.i.d.\ exponential waiting times of mean $\tau = 1$. When a searcher attempts to move, an adjacent node is picked uniformly at random; if the node is vacant the searcher moves there, and if not the move is aborted. Let $T_r^i$ be the time for searcher $i$ to first visit $r$ distinct nodes.
(In practice we evaluate $T_r^i$ by counting the number of attempted jumps made by all parallel searchers until walker~$i$ visits the $r$-th distinct site and rescaling by the average waiting time $1/M$. This is equivalent to describing the CTRW by a discrete-time random walk with time step $1/M$ which is sufficient when seeking only first moments as we do here.)
We define the parallel cover time $\mathcal{C}_M = \text{max}_{1\leq i \leq M} \{T_N^i \}$ as the time for all $M$ searchers to each visit all $N$ nodes, and the expectation $\langle \mathcal{C}_M \rangle$ to be taken simultaneously over the space of all initial configurations and random walk instances. (This is contrary to the usual cover time, where a maximum over initial position is taken~\cite{LovaszRWonGraphs}. Here, this is combinatorially intractable on most networks when $M \gg 1$, and our definition allows us to consider topological effects removed from any influence of initial conditions.)
The APCT, $\mathcal{S}_M = \langle \mathcal{C}_M \rangle / M$, then quantifies the economy of scale in parallel searching: if $\langle \mathcal{C}_M \rangle < M\mathcal{S}_1$ for some $M>1$ there is an efficiency gain in parallel searching beyond the simple serial search strategy. For a given network, let $M^*$ be the optimal number of searchers, that is, the $M$ for which $\mathcal{S}_M$ is minimised, and let $\rho^* = M^*/N$ be the equivalent optimal density. For numerical techniques, see Methods.


\subsection*{Efficient search for ring-like topologies}
Parallel search of a ring lattice is optimal at a strikingly low density ($\rho^*=0.05$ when $N=100$; Fig.~1(a)), becoming increasingly inefficient as $M$ increases. This is due to the strong confinement effects of single-file diffusion~\cite{Beijeren_Tracer3,Lizana_SFDinBox,Lutz_PRL_SFD,Wei_Science_SFD,2017Imamura_PRL,2014Krapivsky_PRL}. That said, parallel searching is still more efficient than simple serial searching for a range of $\rho$ (Fig.~1(a), dashed line). Remarkably, however, a parallel search on a ring lattice can be made significantly more efficient by introducing only a small number of additional edges. To demonstrate, we consider a form of the Newman--Watts ensemble~\cite{NewmanWatts1999} that interpolates between the two extremal topologies of the ring lattice and the complete network. Starting from a ring lattice, we add a fixed number of random additional edges, or shortcuts (Fig.~1(b), inset). The optimal density rapidly approaches that of the complete network, $\rho^*=0.47$, even with only 3--4\% of the possible shortcuts added on $N=100$ nodes (Fig.~1(b)).


Parallel search efficiency on the complete network can be evaluated exactly. In the Supplementary Information (SI Text), we derive the APCT
\begin{equation}\label{eq::int_complete}
\mathcal{S}_M = \dfrac{(N-1)^2}{M(N-M)} h(M(N-1)),
\end{equation}
where $h(i)$ is the $i$-th harmonic number. The optimal density of parallel searchers,~$\rho^*$, for a given $N$ can then be calculated from $M^*=\text{argmin}_M\,\mathcal{S}_M$, which gives ${\rho^* \rightarrow 0.5}$ in the limit $N \rightarrow \infty$ (SI Text).
Equation~\eqref{eq::int_complete} can be written as $\mathcal{S}_M = \phi\tilde{\mathcal{S}}_M$ where $\tilde{\mathcal{S}}_M$ is the non-interacting APCT (SI Text) and $\phi = (1-N^{-1})/(1-\rho)$ is a mean-field correction. This accounts for the average slow-down due to aborted moves by assuming two-site occupancy probabilities to be products of single-site densities (SI Text). While it is exact for the complete network by spatial homogeneity, this need not hold in general.


\begin{figure}[tb]
\centering
\includegraphics[scale=1.3]{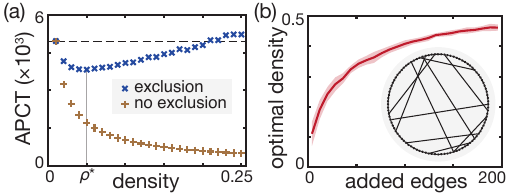}
\caption{(a)~APCT on the ring lattice for interacting ($\times$) and non-interacting ($+$) searchers as a function of searcher density $\rho = M/N$, where $N=100$. The optimal density $\rho^*$ (marked) can be seen as the minimum of the APCT; dashed horizontal line denotes the area below which parallel search is faster than serial search. Each APCT was calculated from $10^5$ random walk instances.
(b)~Mean optimal density of Newman--Watts networks with $N=100$ as a function of the number of shortcuts (example inset). Coloured band indicates $\pm 1$ s.d. Optimal densities were calculated (SI Text) for $1000$ Newman--Watts realisations for a range of added edges.}
\end{figure}


For general networks, na\"{\i}ve mean-field approximations (MFAs)~\cite{LiggetBook,Ghosts_Folks,PlankSimpson_JRSI} can lead to drastically inaccurate estimation of the APCT and optimal density.
Taking the near-Gumbel-distributed single searcher cover times on random networks~\cite{Chupeau_NatPhys_CTs} and attempting to incorporate exclusion through rescaling by $\phi$ fails for many networks: under this MFA the predicted optimal density as $N$ becomes large approaches $\rho^*=0.5$ (SI Text), but for networks of low average degree, and for the ring lattice (Fig.~1(a)) in particular, this is far from the true $\rho^*$ evaluated numerically (Fig.~2).
Insight into this inaccuracy can be drawn from asymptotic estimates of the APCT on the ring lattice in the high density regime. Consider $M=N-k$ searchers with $k \ll N$ vacancies.
Inspired by particle-hole duality~\cite{Beijeren_Tracer3} we can follow the vacancies instead of the searchers. The parallel cover time is then the time taken for the net displacement of the vacancies to first reach $(N-1)(N-k)$~(SI Text). For large $N$ the vacancies are approximately non-interacting, meaning their net displacement is approximately that of a single vacancy moving $k$ times faster, that is, $\tau_{\text{vac}}=1/k$. Standard results~\cite{Yokoi_1990_ExactCTs} then imply
\begin{equation}\label{eq::HighDensity_Symm}
\mathcal{S}_{N-k} \sim \dfrac{(N-1)[(N-1)(N-k)+1]}{2 k}.
\end{equation}
This estimate does extremely well in predicting the APCT (SI Text, Fig.~S2). It is exact for $k=1$, and only noticeably deviates when ${\rho \lesssim 0.85}$ for $N=100$.
Equation~\eqref{eq::HighDensity_Symm} reveals that the APCT at high density is $\mathcal{O}(N^3)$, while the MFA suggests an APCT of $\phi \tilde{\mathcal{S}}_N \sim \mathcal{O}(N^2)$, an order of magnitude difference highlighting the failure of the MFA to capture spatial correlations.


\subsection*{General topologies and the spectral gap}

Can we identify a topological heuristic to predict density optima for general networks? Given the adjacency matrix $\mathbf{A}$ and the diagonal matrix $\mathbf{D}$ whose entries are the node degrees, define the random walk transition matrix $\mathbf{P} = \mathbf{D}^{-1} \mathbf{A}$. The eigenvalues~$\{\lambda_i\}$ of~$\mathbf{P}$ determine how fast the probability distribution of a single-searcher random walk converges to its equilibrium (that is, mixes)~\cite{LovaszRWonGraphs}.
The largest eigenvalue is necessarily $\lambda_1 = 1$, and we define the spectral gap, $\mathcal{G}$, as the difference in magnitude between the first and second largest eigenvalues $\mathcal{G} = 1 - \text{max}_{2 \leqslant i \leqslant N} \{|\lambda_i |\}$.
The closer $\mathcal{G}$ is to zero, the slower a random walk converges~\cite{LovaszRWonGraphs}. The spectral gap is then a natural candidate for predicting density optima as the features of networks that slow down convergence, such as bottlenecks (identified through nodes of high betweenness centrality~\cite{Newman_Networks}), also significantly increase the parallel cover time.
Furthermore, single-searcher cover times can be related~\cite{MaierPRE2017} to the eigenvalues of the combinatorial Laplacian $\mathbf{L} = \mathbf{D} - \mathbf{A}$, which in turn relate to those of ${\mathbf{P} = \mathbf{I} - \mathbf{D}^{-1}\mathbf{L}}$ via the normalised Laplacian~\cite{ChenZhang_07}.


\begin{figure}[tb]
\label{fig:Fig2}
\centering
\includegraphics[scale=1.3]{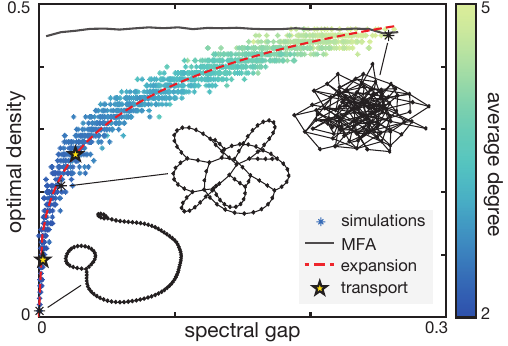}
\caption{Spectral gap predicts optimal parallel search density. Optimal densities were calculated for $1500$ random networks on $N=100$ nodes with average degree between $2$ and $5$ (Methods). Each point represents one network with average degree indicated by colour, with three topologically contrasting networks highlighted. Solid grey curve denotes the mean MFA-predicted optimal density as a function of $\mathcal{G}$ from numerically-determined single-searcher cover times. Dashed red curve shows best-fit expansion $\rho^*(\mathcal{G};0.7556,0.2933)$. Two real-world transport networks are indicated by stars: the London Underground (lower) and American airports (upper).}
\end{figure}


There is a tight relationship between the spectral gap and the optimal density of searchers for general networks (Fig.~2),  quantified through high mutual information (SI Text). Sampling over random networks of minimum degree~$2$ and average degree between $2$ and $5$ (Methods) reveals that topology has a huge impact on optimal parallel search strategies: for a small spectral gap (Fig.~2, right-most inset) networks have low average degree and a high concentration of nodes of degree $2$, resulting in linear chains along which the single-file diffusion of searchers significantly increases the cover time \cite{Beijeren_Tracer3,Lizana_SFDinBox,Lutz_PRL_SFD,Wei_Science_SFD,2017Imamura_PRL,2014Krapivsky_PRL}.
Indeed, for a linear chain the parallel cover time is finite only for a single searcher; for two or more searchers the left-most (right-most)  searcher will never reach the right-most (left-most)  site due to the effect of single-file diffusion.
The na\"{\i}ve MFA fails to capture this relationship (Fig.~2, grey curve), demonstrating the importance of long-lived occupancy correlations.
The MFA only becomes valid for networks of average degree nearing $5$, where the density optima approach that of the complete network ($\rho^*=0.47$ for $N=100$). Typically, these networks have a small fraction of degree two nodes and are more highly connected (Fig.~2, left-most inset). A power series $\rho^*\left( \mathcal{G};a,b\right) = \mathcal{G}^b [ 1/2 - \left(1/2 - a \right) \left(1-\mathcal{G}\right)]$, constructed such that $\rho^*(1) = 1/2$ to match the $N \rightarrow \infty$ limit of the complete graph, matches the data well (Fig.~2, red curve) and provides an easily-computed predictive approximation.
Beyond random networks, we find that two real-world transport networks, the London Underground and American airport networks~\cite{NatPhysAirport}, further validate this relationship in practice (Fig.~2). While the inefficiency of the former may be expected from an average degree near~$2$, the latter is more surprising: despite an average degree near~$12$ its optimal density is markedly low, a phenomenon captured by the spectral gap~(SI Text).
Simpler statistics based on the degree distribution fail in such cases and have lower mutual information over our graph ensemble (SI Text).


\begin{figure}[tb]
\label{fig:Fig3}
\centering
\includegraphics[scale=1.3]{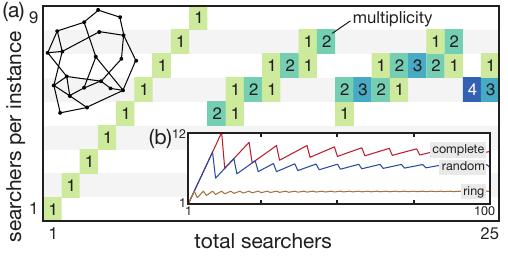}
\caption{(a)~Optimal search strategy as $M_{\text{tot}}$ is increased for a random network of $N= 20$ nodes (inset) with $\rho^*=0.35$. Values at fixed $M_\text{tot}$ give the non-zero components of the  strategy vector $\mathbf{k}_\text{opt}$. (b)~Mean number of searchers per instance $\langle \mathbf{k}_\text{opt} \rangle$ for the optimal search strategies of the complete network, the random network in (a) and the ring lattice on $N=20$ nodes, showing the approach of $\langle \mathbf{k}_\text{opt} \rangle$ towards the optimal APCT densities $\rho^* = 0.45$, $0.35$ and $0.15$ respectively for large $M_\text{tot}$. APCTs are averages over $10^5$ random walk instances.}
\end{figure}

\subsection*{Optimising parallel search strategies}

Our results enable the design of optimal search strategies for arbitrary numbers of searchers.
Suppose $M_{\text{tot}}$ searchers, potentially more than $N$, need to each cover a given network.
A hybrid series--parallel strategy comprises sequentially introducing batches of searchers, with each batch covering the network in parallel to completion before being removed and replaced by the next batch.
Let $k_i$ be the number of times~$i$ parallel searchers are introduced sequentially; for example, $k_1=3$ and $k_2=4$ denotes performing three successive searches with one searcher and four successive searches with two searchers. A hybrid strategy is then defined by the vector $\mathbf{k} = \left( k_1,\ldots,k_{N-1}\right)$. The optimal strategy, $\mathbf{k}_\text{opt}$, minimises the overall cover time $\mathcal{T} =\sum_{i=1}^{N-1} k_i \langle \mathcal{C}_i \rangle$
subject to the constraint $\sum_{i=1}^{N-1} i k_i = M_{\text{tot}}$. As $M_{\text{tot}}$ increases the solution of this integer programming problem typically yields a strategy comprising a mixture of searcher numbers near the optimal density (Fig.~3). For $M_{\text{tot}} \gg N$, we find that the optimal strategy is to have almost all serial searches operating at optimal density with vanishingly small standard deviation away from the optimum (Fig.~3(b); SI Text), validating optimal density as a crucial search statistic.
Thus for large networks whose optimal density becomes difficult to evaluate, our results present a simple strategy to optimise mutually excluding search: approximate $\rho^*$ via the easily-computed spectral gap (Fig.~2), and then perform searches in batches at this density.

\subsection*{Asymmetric search processses}

Until now, we have only considered the symmetric exclusion process, where searchers uniformly sample a target node from their neighbours. However, many biological and physical systems are not in equilibrium, possessing non-zero probability currents in stationary state that markedly change the possible physical behaviours~\cite{2011Neri_PRL}. To model this, we now explore directed networks with bias in the choice of target node. To avoid searchers accumulating at nodes we restrict to flux-conserving (balanced) networks, where each vertex has an equal number of inwards- and outwards-biased edges (and is therefore necessarily of even degree; Fig.~4(a)). For each edge biased $u \rightarrow v$, we define the transition probabilities as $p_{u\rightarrow v} = \left(1 + 2\varepsilon \right)/d_u$ and $p_{v \rightarrow u} = \left( 1 - 2\varepsilon \right)/d_v$, where $d_u$, $d_v$ are the vertex degrees and $\varepsilon \in [0,1/2]$ controls the bias. Thus the probability a searcher exits a node through an outwardly biased edge is $1/2 + \varepsilon$.


Asymmetric processes facilitate a significantly greater optimal density than the equivalent unbiased process on the ring lattice.
We find numerically that the optimal density for $N=100$ and $\varepsilon = 0.25$ is $\rho^* = 0.51$ (with high confidence; SI Text) compared to $\rho^*=0.05$ in the symmetric case $\varepsilon=0$. This is a marked increase in $\rho^*$ that even exceeds the symmetrically-biased complete network for which $\rho^* < 0.5$. Here, conservation of flux forces there to be only two viable orientations for the edges, all either clockwise or anticlockwise, effectively adding a constant drift to the symmetric process.
As for symmetric exclusion, intuition can be drawn from asymptotic estimates of the APCT in the high density regime for the ring lattice. Let $p$ and $q = 1-p$ be the probabilities of moving clockwise and anticlockwise, respectively. In the SI Text we show that, for large $N$ and small vacancy count $k$, the APCT for $M=N-k$ searchers with $p>q$ is approximately $ \mathcal{S}_{N-k} \sim (N-1)/(p-q)k$. This is $\mathcal{O}(N)$ for a non-negligible bias, two orders of magnitude smaller than Eq.~\eqref{eq::HighDensity_Symm}. In addition, the APCTs in the low and high density regimes are both $\mathcal{O}(N)$ (SI Text) in contrast to the unbiased process, whose low and high density APCTs we recall as $\mathcal{O}(N^2)$ and $\mathcal{O}(N^3)$, respectively, indicating bias induces a qualitative change in the APCT.


\begin{figure}[tb]
\label{fig:Fig4}
\centering
\includegraphics[scale=1.3]{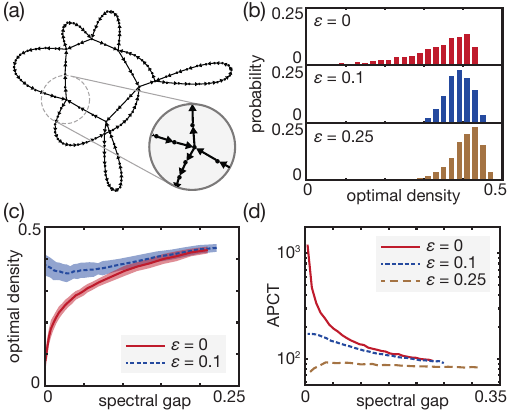}
\caption{
(a)~A random network for the asymmetric exclusion process, where $N=100$. The arrows depict the directional bias of the random searchers between nodes.
(b)~Frequency of optimal densities for a variety of biases.
(c)~Mean optimal density as a function of spectral gap for symmetric (${\varepsilon=0}$) and asymmetric (${\varepsilon=0.1}$) exclusion processes. Coloured bands are $\pm 1$ s.d.
(d)~APCT at optimal searcher density for a variety of biases.
Results in (b) and (d) are from the same ensemble of $2250$ random networks of average degree between $2$ and $5$ (Methods). In (c) an additional $2000$ networks are included with average degree between $2$ and $3$. APCTs in (d) are averages over $10^5$ random walk instances.}
\end{figure}


Beyond the ring lattice, general flux-conserving networks of low degree also see a significantly enhanced optimal density for sufficient bias, with a far narrower distribution of optimal densities (Fig.~4(b)). Dramatic increases in optimal density particularly occur for networks with $\mathcal{G}\approx 0$ (Fig.~4(c)): these have a high concentration of degree two nodes implying many linear chains, the edges of which must have the same directional bias (Fig.~4(a)). More surprisingly, for networks with extremely small spectral gaps ($\mathcal{G} \approx 0$) we see a non-monotonic relationship between $\mathcal{G}$ and $\rho^*$ on average. This phenomenon is also captured by the APCT at optimal density (Fig.~4(d)) where networks with  $\mathcal{G} \approx 0$ typically have APCTs below those of networks with greater~$\mathcal{G}$. This shows that the addition of flux-conserving bias can have a counterintuitive impact on search efficiency and optimality (which can be orientation dependent; SI Text).
For example, consider $M_\text{tot} = 100$ searchers on the random network in Fig.~4(a) using the optimal serial--parallel strategy versus the same search on the ring lattice. Without bias ($\varepsilon = 0$), the random network has an average search time of $6.6 \times 10^4$ and the ring lattice has time $4.8 \times 10^5$. In contrast, with $\varepsilon = 0.25$, these times become $8.5 \times 10^3$ and $8.3 \times 10^2$, respectively. Thus, random flux-conserving bias inverts the networks' relative efficiency, implying  care must be taken when attempting to improve network search efficiency process by na\"{\i}ve biasing.

\section*{Conclusion}

We have introduced the APCT as a fundamental measure of how efficiently a network can support mutually excluding random search processes. We found that the optimal density of concurrent searchers can vary from $1$--$2\%$ to nearly $50\%$ depending on topology, and showed that this topological dependence can be efficiently captured by the spectral gap for both artificial random graphs and real-world transport networks.
Finally, we generalised to a biased, non-equilibrium process and uncovered a qualitative change in the topological dependence of optimal densities. This evidences how the transport process itself must be taken into account when optimising random search, and leads us to ask whether there is an efficient  statistical characterisation of the topology--transport interplay for general active networked processes~\cite{2000Caspi_PRL,2013Becker_PRL,2016Woodhouse_PNAS,2017Graf_PRL}. More broadly, our work paves the way towards new strategies for topology optimisation in process allocation and flow transport problems across physical and biological networked systems.


\section*{Methods}

\subsection*{Numerical evaluation of cover time statistics and optimal densities}
All cover time statistics presented in this work were calculated as follows. We do not explicitly sample from the exponential waiting time distribution; rather, we sample the next attempted jump and update the time statistic of interest, $t$, by the average waiting time, $t = t + \tau / M$. As we are only interested in the first moments of the cover times, linearity of expectation implies that this suffices. We now detail our numerical procedure for evaluating optimal densities. Suppose we have a network with $N$ nodes and we wish to calculate the optimal number of parallel searchers, $M^*$, that minimises the APCT. For efficiency we begin using a small number of random walk instances (initially $10^3$) that we increase as we get closer to $M^*$. Starting at $M=1$ we calculate the average cover time $\mathcal{S}_1$, then increase $M$ by one and run another sample of random walk instances, calculating $\mathcal{S}_2$ and so on. After generating a new APCT for $j+1$ parallel searchers we ask whether $\mathcal{S}_{j+1} > \mathcal{S}_j$ to see whether $j$ is a candidate for $M^*$. If not we increase $M$ by one, however if so then we increase the number of random walk instances to $10^4$ and return to the previous number of parallel searchers $j$, where we recalculate $\mathcal{S}_j$ and $\mathcal{S}_{j+1}$, and check again to see if $\mathcal{S}_{j+1} > \mathcal{S}_j$. Once we have moved to $10^4$ random walk instances we still check to see if the new APCT is greater than the old one, as before. Once we find a density of parallel searchers that appears to be the minimum we increase the number of random walk instances used for the final time to $10^5$. Then the number of searchers that are found to minimise the APCT using $10^5$ random walk instances is taken to be $M^*$, the optimal number of searchers. This was typically enough realisations such that the standard error when calculating the APCTs resulted in mutually excluding 95\% confidence intervals either side of the optimal density.

\subsection*{Generating random symmetric networks}

In this section we explain the procedure used to sample the random networks in Fig.~2 of the text.

\begin{algorithmic}[1]

  \State{Select $m$ and $N$, the required number of edges and the number of nodes of the network respectively, such that $m \geq N$.}

  \State{Sample uniformly the target degree distribution, $\vec{d}$, such that every node has degree at least two, i.e. $\vec{d} = 2(1,\ldots,1) + \mathrm{Multinomial}(2(m-N),N)$.}
\label{degree_dist}

  \State{Allocate $d_i$ stubs to the $i$-th node for $i\in \{1,\ldots,N\}$.}

  \While{ there are stubs remaining }

  \If{the only remaining connections to choose from result in self-loops or multiple-edges } start over and return to step \ref{degree_dist}.

  \Else{ select two distinct nodes (avoiding self-loops) with stubs remaining uniformly at random and connect them with an edge only if they have not already been connected (avoiding multiple-edges).}

  \EndIf

  \EndWhile

  \If{the resulting network is not connected} return to step \ref{degree_dist}.
  \Else{  take the resulting network to be a single realisation.}
  \EndIf
\end{algorithmic}
\color{black}

\subsection*{Generating random balanced directed networks}

In this section we describe the sampling procedure for the asymmetric random networks in Fig.~4 of the text.

\begin{algorithmic}[1]

  \State{Select $m$ and $N$, the required number of edges and the number of nodes of the network respectively, such that $m \geq N$.}

  \State{Sample uniformly the target degree distribution for outwardly biased edges, $\vec{d}_{\text{out}}$, such that every node has outward degree at least one, i.e. $\vec{d}_{\text{out}} = (1,\ldots,1) + \mathrm{Multinomial}(m-N,N)$.}
\label{degree_dist_asym}

\State{Set the degree distribution for the inwardly biased edges to be $\vec{d}_{\text{in}}=\vec{d}_{\text{out}}$.}

  \State{Allocate $\vec{d}_{\text{out}}(i)=\vec{d}_{\text{in}}(i)$ outward stubs and inward stubs to the $i$-th node for $i\in \{1,\ldots,N\}$.}

  \While{ there are stubs remaining }

  \If{the only remaining connections to choose from result in self-loops or multiple-edges } start over and return to step \ref{degree_dist_asym}.

  \Else{ select uniformly at random a node with an outward stub and a distinct node with an inward stub. Connect these two nodes with an edge with preferential direction from the outward node to the inward node as long as these nodes have not been connected before.}

  \EndIf

  \EndWhile

  \If{the resulting network is not connected} return to step \ref{degree_dist_asym}.
  \Else{ take the resulting network to be a single realisation.}
  \EndIf
\end{algorithmic}
\color{black}


\section*{Acknowledgments}
This work was supported by the EPSRC Systems Biology DTC EP/G03706X/1 (D.B.W.), a Royal Society Wolfson Research Merit Award (R.E.B.), a Leverhulme Research Fellowship (R.E.B.), the BBSRC Multi-Scale Biology Network (R.E.B. and F.G.W.), and Trinity College, Cambridge (F.G.W.).


\section*{Competing Interests}

The authors declare no competing financial interests.


\section*{Author Contributions}

All authors contributed at all stages of this work.


\singlespacing


\begin{thebibliography}{10}
\expandafter\ifx\csname url\endcsname\relax
  \def\url#1{\texttt{#1}}\fi
\expandafter\ifx\csname urlprefix\endcsname\relax\def\urlprefix{URL }\fi
\providecommand{\bibinfo}[2]{#2}
\providecommand{\eprint}[2][]{\url{#2}}

\bibitem{Benichou_ISS_review}
\bibinfo{author}{B\'enichou, O.}, \bibinfo{author}{Loverdo, C.},
  \bibinfo{author}{Moreau, M.} \& \bibinfo{author}{Voituriez, R.}
\newblock \bibinfo{title}{Intermittent search strategies}.
\newblock \emph{\bibinfo{journal}{Rev. Mod. Phys.}}
  \textbf{\bibinfo{volume}{83}}, \bibinfo{pages}{81--129}
  (\bibinfo{year}{2011}).

\bibitem{Optimal_RunandTumble}
\bibinfo{author}{Rupprecht, J.-F.}, \bibinfo{author}{B\'enichou, O.} \&
  \bibinfo{author}{Voituriez, R.}
\newblock \bibinfo{title}{Optimal search strategies of run-and-tumble walks}.
\newblock \emph{\bibinfo{journal}{Phys. Rev. E}} \textbf{\bibinfo{volume}{94}},
  \bibinfo{pages}{012117} (\bibinfo{year}{2016}).

\bibitem{Coppey_ProteinDNA}
\bibinfo{author}{Coppey, M.}, \bibinfo{author}{B\'enichou, O.},
  \bibinfo{author}{Voituriez, R.} \& \bibinfo{author}{Moreau, M.}
\newblock \bibinfo{title}{Kinetics of target site localization of a protein on
  dna: A stochastic approach}.
\newblock \emph{\bibinfo{journal}{Biophys. J.}} \textbf{\bibinfo{volume}{87}},
  \bibinfo{pages}{1640--1649} (\bibinfo{year}{2004}).

\bibitem{Bartumeus_2009}
\bibinfo{author}{Bartmeus, F.} \& \bibinfo{author}{Catalan, J.}
\newblock \bibinfo{title}{Optimal search behaviour and classical foraging
  theory}.
\newblock \emph{\bibinfo{journal}{J. Phys. A-Math. Theor.}}
  \textbf{\bibinfo{volume}{42}}, \bibinfo{pages}{434002}
  (\bibinfo{year}{2009}).

\bibitem{pathogens}
\bibinfo{author}{Heuz\'e, M.~L.} \emph{et~al.}
\newblock \bibinfo{title}{Migration of dendritic cells: Physical principles,
  molecular mechanisms, and function implications.}
\newblock \emph{\bibinfo{journal}{Immunol. Rev.}}
  \textbf{\bibinfo{volume}{256}}, \bibinfo{pages}{240} (\bibinfo{year}{2013}).

\bibitem{Levy_foraging}
\bibinfo{author}{Viswanathan, G.~M.}, \bibinfo{author}{Raposo, E.~P.} \&
  \bibinfo{author}{da~Luz, M. G.~E.}
\newblock \bibinfo{title}{L\'evy flights and superdiffusion in the context of
  biological encounters and random searches}.
\newblock \emph{\bibinfo{journal}{Phys. Life Rev.}}
  \textbf{\bibinfo{volume}{5}}, \bibinfo{pages}{133} (\bibinfo{year}{2008}).

\bibitem{OptimalRWs_Nature}
\bibinfo{author}{Viswanathan, G.~M.} \emph{et~al.}
\newblock \bibinfo{title}{Optimizing the success of random searches}.
\newblock \emph{\bibinfo{journal}{Nature}} \textbf{\bibinfo{volume}{401}},
  \bibinfo{pages}{911} (\bibinfo{year}{1999}).

\bibitem{RW_in_Bio_JRSI}
\bibinfo{author}{Codling, E.~A.}, \bibinfo{author}{Plank, M.~J.} \&
  \bibinfo{author}{Benhamou, S.}
\newblock \bibinfo{title}{Random walks in biology}.
\newblock \emph{\bibinfo{journal}{J. R. Soc. Interface}}
  \textbf{\bibinfo{volume}{5}}, \bibinfo{pages}{813} (\bibinfo{year}{2008}).

\bibitem{2009Li_NatPhys}
\bibinfo{author}{Li, G.-W.}, \bibinfo{author}{Berg, O.~G.} \&
  \bibinfo{author}{Elf, J.}
\newblock \bibinfo{title}{Effects of macromolecular crowding and {DNA} looping
  on gene regulation kinetics}.
\newblock \emph{\bibinfo{journal}{Nat. Phys.}} \textbf{\bibinfo{volume}{5}},
  \bibinfo{pages}{294} (\bibinfo{year}{2009}).

\bibitem{2016Schwarz_PRL}
\bibinfo{author}{Schwarz, K.}, \bibinfo{author}{Schr{\"o}der, Y.},
  \bibinfo{author}{Qu, B.}, \bibinfo{author}{Hoth, M.} \&
  \bibinfo{author}{Rieger, H.}
\newblock \bibinfo{title}{Optimality of spatially inhomogeneous search
  strategies}.
\newblock \emph{\bibinfo{journal}{Phys. Rev. Lett.}}
  \textbf{\bibinfo{volume}{117}}, \bibinfo{pages}{068101}
  (\bibinfo{year}{2016}).

\bibitem{2014Rubenstein_Science}
\bibinfo{author}{Rubenstein, M.}, \bibinfo{author}{Cornejo, A.} \&
  \bibinfo{author}{Nagpal, R.}
\newblock \bibinfo{title}{Programmable self-assembly in a thousand-robot
  swarm}.
\newblock \emph{\bibinfo{journal}{Science}} \textbf{\bibinfo{volume}{345}},
  \bibinfo{pages}{795--799} (\bibinfo{year}{2014}).

\bibitem{Herty_TrafficFlow}
\bibinfo{author}{Herty, M.} \& \bibinfo{author}{Klar, A.}
\newblock \bibinfo{title}{Modelling, simulation, and optimization of traffic
  flow networks}.
\newblock \emph{\bibinfo{journal}{SIAM J. Sci. Comput.}}
  \textbf{\bibinfo{volume}{25}}, \bibinfo{pages}{1066} (\bibinfo{year}{2008}).

\bibitem{SIAMRev_Bellomo}
\bibinfo{author}{Bellomo, N.} \& \bibinfo{author}{Dogbe, C.}
\newblock \bibinfo{title}{On the modelling of traffic and crowds: a survey of
  models, speculations, and perspectives}.
\newblock \emph{\bibinfo{journal}{SIAM Rev.}} \textbf{\bibinfo{volume}{53}},
  \bibinfo{pages}{409} (\bibinfo{year}{2011}).

\bibitem{SIAMFugen}
\bibinfo{author}{F\"ugenschuh, A.}, \bibinfo{author}{Herty, M.},
  \bibinfo{author}{Klar, A.} \& \bibinfo{author}{Martin, A.}
\newblock \bibinfo{title}{Combinatorial and continuous models for the
  optimization of traffic flows on networks}.
\newblock \emph{\bibinfo{journal}{SIAM J. Optim.}}
  \textbf{\bibinfo{volume}{16}}, \bibinfo{pages}{1155} (\bibinfo{year}{2006}).

\bibitem{2001Adamic_PRE}
\bibinfo{author}{Adamic, L.~A.}, \bibinfo{author}{Lukose, R.~M.},
  \bibinfo{author}{Puniyani, A.~R.} \& \bibinfo{author}{Huberman, B.~A.}
\newblock \bibinfo{title}{Search in power-law networks}.
\newblock \emph{\bibinfo{journal}{Phys. Rev. E}} \textbf{\bibinfo{volume}{64}},
  \bibinfo{pages}{046135} (\bibinfo{year}{2001}).

\bibitem{Chupeau_NatPhys_CTs}
\bibinfo{author}{Chupeau, M.}, \bibinfo{author}{B\'enichou, O.} \&
  \bibinfo{author}{Voituriez, R.}
\newblock \bibinfo{title}{Cover times of random searches}.
\newblock \emph{\bibinfo{journal}{Nat. Phys.}} \textbf{\bibinfo{volume}{11}},
  \bibinfo{pages}{844--847} (\bibinfo{year}{2015}).

\bibitem{Kahn1989_CTs}
\bibinfo{author}{Kahn, J.~D.}, \bibinfo{author}{Linial, N.},
  \bibinfo{author}{Nisan, N.} \& \bibinfo{author}{Saks, M.~E.}
\newblock \bibinfo{title}{On the cover time of random walks on graphs.}
\newblock \emph{\bibinfo{journal}{J. Theor. Prob.}}
  \textbf{\bibinfo{volume}{2}}, \bibinfo{pages}{121} (\bibinfo{year}{1989}).

\bibitem{Kundu_ExactDistinctCommon}
\bibinfo{author}{Kundu, A.}, \bibinfo{author}{Majumdar, S.~N.} \&
  \bibinfo{author}{Schehr, G.}
\newblock \bibinfo{title}{Exact distributions of the number of distinct and
  common sites visited by n independent random walkers}.
\newblock \emph{\bibinfo{journal}{Phys. Rev. Lett.}}
  \textbf{\bibinfo{volume}{110}}, \bibinfo{pages}{220602}
  (\bibinfo{year}{2013}).

\bibitem{Majumdar_ExactCTRing}
\bibinfo{author}{Majumdar, S.~N.}, \bibinfo{author}{Sabhapandit, S.} \&
  \bibinfo{author}{Schehr, G.}
\newblock \bibinfo{title}{Exact distributions of cover times or n independent
  random walkers in one dimension}.
\newblock \emph{\bibinfo{journal}{Phys. Rev. E}} \textbf{\bibinfo{volume}{94}},
  \bibinfo{pages}{062131} (\bibinfo{year}{2016}).

\bibitem{Yokoi_1990_ExactCTs}
\bibinfo{author}{Yokoi, C. S.~O.}, \bibinfo{author}{Hern\'andez-Machado, A.} \&
  \bibinfo{author}{Ram\'irez-Piscina, L.}
\newblock \bibinfo{title}{Some exact results for the lattice covering time
  problem}.
\newblock \emph{\bibinfo{journal}{Phys. Lett. A}}
  \textbf{\bibinfo{volume}{145}}, \bibinfo{pages}{82--86}
  (\bibinfo{year}{1990}).

\bibitem{CTs_2Dlattice}
\bibinfo{author}{Ding, J.}
\newblock \bibinfo{title}{On cover times for 2d lattices}.
\newblock \emph{\bibinfo{journal}{Electron. J. Probab.}}
  \textbf{\bibinfo{volume}{17}}, \bibinfo{pages}{1} (\bibinfo{year}{2012}).

\bibitem{CTs_Lattice}
\bibinfo{author}{Brummelhuis, M. J. A.~M.} \& \bibinfo{author}{Hilhorst, H.~J.}
\newblock \bibinfo{title}{Covering of a finite lattice by a random walk}.
\newblock \emph{\bibinfo{journal}{Physica A}} \textbf{\bibinfo{volume}{176}},
  \bibinfo{pages}{387} (\bibinfo{year}{1991}).

\bibitem{MaierPRE2017}
\bibinfo{author}{Maier, B.~F.} \& \bibinfo{author}{Brockmann, D.}
\newblock \bibinfo{title}{Cover time for random walks on arbitrary complex
  networks}.
\newblock \emph{\bibinfo{journal}{Phys. Rev. E}} \textbf{\bibinfo{volume}{96}},
  \bibinfo{pages}{042307} (\bibinfo{year}{2017}).

\bibitem{Tejedor_PRE09}
\bibinfo{author}{Tejedor, V.}, \bibinfo{author}{B\'enichou, O.} \&
  \bibinfo{author}{Voituriez, R.}
\newblock \bibinfo{title}{Global mean first-passage times of random walks on
  complex networks}.
\newblock \emph{\bibinfo{journal}{Phys. Rev. E}} \textbf{\bibinfo{volume}{80}},
  \bibinfo{pages}{065194} (\bibinfo{year}{2009}).

\bibitem{Meyer_PRE12}
\bibinfo{author}{Meyer, B.}, \bibinfo{author}{Agliari, E.},
  \bibinfo{author}{B\'enichou, O.} \& \bibinfo{author}{Voituriez, R.}
\newblock \bibinfo{title}{Exact calculations of first-passage quantites on
  recursive networks}.
\newblock \emph{\bibinfo{journal}{Phys. Rev. E}} \textbf{\bibinfo{volume}{85}},
  \bibinfo{pages}{026113} (\bibinfo{year}{2012}).

\bibitem{Palacios_03}
\bibinfo{author}{Palacios, J.}
\newblock \bibinfo{title}{On hitting and cover times via electric networks}.
\newblock \emph{\bibinfo{journal}{Braz. J. Probab. Stat.}}
  \textbf{\bibinfo{volume}{17}}, \bibinfo{pages}{127} (\bibinfo{year}{2003}).

\bibitem{JonassonCTs}
\bibinfo{author}{Jonasson, J.}
\newblock \bibinfo{title}{On the cover times for random walks on random
  graphs}.
\newblock \emph{\bibinfo{journal}{Comb. Probab. Comput.}}
  \textbf{\bibinfo{volume}{7}}, \bibinfo{pages}{265} (\bibinfo{year}{1998}).

\bibitem{2017Weng_PRE}
\bibinfo{author}{Weng, T.}, \bibinfo{author}{Zhang, J.},
  \bibinfo{author}{Small, M.} \& \bibinfo{author}{Hui, P.}
\newblock \bibinfo{title}{Multiple random walks on complex networks: {A}
  harmonic law predicts search time}.
\newblock \emph{\bibinfo{journal}{Phys. Rev. E}} \textbf{\bibinfo{volume}{95}},
  \bibinfo{pages}{052103} (\bibinfo{year}{2017}).

\bibitem{LiggetBook}
\bibinfo{author}{Liggett, T.}
\newblock \emph{\bibinfo{title}{Stochastic Interacting Systems: Contact, Voter
  and Exclusion Processes}} (\bibinfo{publisher}{Springer},
  \bibinfo{year}{1999}).

\bibitem{Ghosts_Folks}
\bibinfo{author}{Simpson, M.~J.}, \bibinfo{author}{Hughes, B.~D.} \&
  \bibinfo{author}{Landman, K.~A.}
\newblock \bibinfo{title}{Diffusing populations: ghosts or folks?}
\newblock \emph{\bibinfo{journal}{Australas. J. Eng. Ed.}}
  \textbf{\bibinfo{volume}{15}}, \bibinfo{pages}{59} (\bibinfo{year}{2009}).

\bibitem{PlankSimpson_JRSI}
\bibinfo{author}{Plank, M.~J.} \& \bibinfo{author}{Simpson, M.~J.}
\newblock \bibinfo{title}{Models of collective cell behaviour with crowding
  effects: comparing lattice-based and lattice-free approaches}.
\newblock \emph{\bibinfo{journal}{J. R. Soc. Interface}}
  \textbf{\bibinfo{volume}{9}}, \bibinfo{pages}{2983} (\bibinfo{year}{2012}).

\bibitem{2005Chowdhury_PhysLifeRev}
\bibinfo{author}{Chowdhury, D.}, \bibinfo{author}{Schadschneider, A.} \&
  \bibinfo{author}{Nishinari, K.}
\newblock \bibinfo{title}{Physics of transport and traffic phenomena in
  biology: from molecular motors and cells to organisms}.
\newblock \emph{\bibinfo{journal}{Phys. Life Rev.}}
  \textbf{\bibinfo{volume}{2}}, \bibinfo{pages}{318--352}
  (\bibinfo{year}{2005}).

\bibitem{2017Graf_PRL}
\bibinfo{author}{Graf, I.~R.} \& \bibinfo{author}{Frey, E.}
\newblock \bibinfo{title}{Generic transport mechanisms for molecular traffic in
  cellular protrusions}.
\newblock \emph{\bibinfo{journal}{Phys. Rev. Lett.}}
  \textbf{\bibinfo{volume}{118}}, \bibinfo{pages}{128101}
  (\bibinfo{year}{2017}).

\bibitem{2011Neri_PRL}
\bibinfo{author}{Neri, I.}, \bibinfo{author}{Kern, N.} \&
  \bibinfo{author}{Parmeggiani, A.}
\newblock \bibinfo{title}{Totally asymmetric simple exclusion process on
  networks}.
\newblock \emph{\bibinfo{journal}{Phys. Rev. Lett.}}
  \textbf{\bibinfo{volume}{107}}, \bibinfo{pages}{068702}
  (\bibinfo{year}{2011}).

\bibitem{1987Plischke_PRB}
\bibinfo{author}{Plischke, M.}, \bibinfo{author}{R{\'a}cz, Z.} \&
  \bibinfo{author}{Liu, D.}
\newblock \bibinfo{title}{Time-reversal invariance and universality of
  two-dimensional growth models}.
\newblock \emph{\bibinfo{journal}{Phys. Rev. B}} \textbf{\bibinfo{volume}{35}},
  \bibinfo{pages}{3485} (\bibinfo{year}{1987}).

\bibitem{1986Kipnis_AnnProb}
\bibinfo{author}{Kipnis, C.}
\newblock \bibinfo{title}{Central limit theorems for infinite series of queues
  and applications to simple exclusion}.
\newblock \emph{\bibinfo{journal}{Ann. Probab.}} \textbf{\bibinfo{volume}{14}},
  \bibinfo{pages}{397--408} (\bibinfo{year}{1986}).

\bibitem{LovaszRWonGraphs}
\bibinfo{author}{Lov\'asz, L.}
\newblock \bibinfo{title}{Random walks on graphs: A survey}.
\newblock \emph{\bibinfo{journal}{Combinatorics, Paul Erd\"os is eighty}}
  \textbf{\bibinfo{volume}{2}}, \bibinfo{pages}{1--46} (\bibinfo{year}{1993}).

\bibitem{Beijeren_Tracer3}
\bibinfo{author}{van Beijeren, H.}, \bibinfo{author}{Kehr, K.~W.} \&
  \bibinfo{author}{Kutner, R.}
\newblock \bibinfo{title}{Diffusion in concentrated lattice gases. iii. tracer
  diffusion on a one-dimensional lattice}.
\newblock \emph{\bibinfo{journal}{Phys. Rev. B}} \textbf{\bibinfo{volume}{28}},
  \bibinfo{pages}{5811} (\bibinfo{year}{1983}).

\bibitem{Lizana_SFDinBox}
\bibinfo{author}{Lizana, L.} \& \bibinfo{author}{Ambj\"{o}rnsson, T.}
\newblock \bibinfo{title}{Single-file diffusion in a box}.
\newblock \emph{\bibinfo{journal}{Phys. Rev. Lett.}}
  \textbf{\bibinfo{volume}{100}}, \bibinfo{pages}{200601}
  (\bibinfo{year}{2008}).

\bibitem{Lutz_PRL_SFD}
\bibinfo{author}{Lutz, C.}, \bibinfo{author}{Kollmann, M.} \&
  \bibinfo{author}{Bechinger, C.}
\newblock \bibinfo{title}{Single-file diffusion of colloids in one-dimensional
  channels}.
\newblock \emph{\bibinfo{journal}{Phys. Rev. Lett.}}
  \textbf{\bibinfo{volume}{93}}, \bibinfo{pages}{026001}
  (\bibinfo{year}{2004}).

\bibitem{Wei_Science_SFD}
\bibinfo{author}{Wei, Q.}, \bibinfo{author}{Bechinger, C.} \&
  \bibinfo{author}{Leiderer, P.}
\newblock \bibinfo{title}{Single-file diffusion of colloids in one-dimensional
  channels}.
\newblock \emph{\bibinfo{journal}{Science}} \textbf{\bibinfo{volume}{287}},
  \bibinfo{pages}{625} (\bibinfo{year}{2000}).

\bibitem{2017Imamura_PRL}
\bibinfo{author}{Imamura, T.}, \bibinfo{author}{Mallick, K.} \&
  \bibinfo{author}{Sasamoto, T.}
\newblock \bibinfo{title}{Large deviations of a tracer in the symmetric
  exclusion process}.
\newblock \emph{\bibinfo{journal}{Phys. Rev. Lett.}}
  \textbf{\bibinfo{volume}{118}}, \bibinfo{pages}{160601}
  (\bibinfo{year}{2017}).

\bibitem{2014Krapivsky_PRL}
\bibinfo{author}{Krapivsky, P.~L.}, \bibinfo{author}{Mallick, K.} \&
  \bibinfo{author}{Sadhu, T.}
\newblock \bibinfo{title}{Large deviations in single-file diffusion}.
\newblock \emph{\bibinfo{journal}{Phys. Rev. Lett.}}
  \textbf{\bibinfo{volume}{113}}, \bibinfo{pages}{078101}
  (\bibinfo{year}{2014}).

\bibitem{NewmanWatts1999}
\bibinfo{author}{Newman, M. E.~J.} \& \bibinfo{author}{Watts, D.~J.}
\newblock \bibinfo{title}{Renormalization group analysis of the small-world
  network model}.
\newblock \emph{\bibinfo{journal}{Phys. Lett. A}}
  \textbf{\bibinfo{volume}{263}}, \bibinfo{pages}{341} (\bibinfo{year}{1999}).

\bibitem{Newman_Networks}
\bibinfo{author}{Newman, M. E.~J.}
\newblock \emph{\bibinfo{title}{Networks: An Introduction}}
  (\bibinfo{publisher}{Oxford University Press}, \bibinfo{year}{2010}).

\bibitem{ChenZhang_07}
\bibinfo{author}{Chen, H.} \& \bibinfo{author}{Zhang, F.}
\newblock \bibinfo{title}{Resistance distance and the normalized {L}aplacian
  spectrum}.
\newblock \emph{\bibinfo{journal}{Discrete Appl. Math.}}
  \textbf{\bibinfo{volume}{155}}, \bibinfo{pages}{654--661}
  (\bibinfo{year}{2007}).

\bibitem{NatPhysAirport}
\bibinfo{author}{Colizza, V.}, \bibinfo{author}{Pastor-Sattoras, R.} \&
  \bibinfo{author}{Vespignani, A.}
\newblock \bibinfo{title}{Reaction--diffusion processes and metapopulation
  models in heterogeneous networks}.
\newblock \emph{\bibinfo{journal}{Nat. Phys.}} \textbf{\bibinfo{volume}{3}},
  \bibinfo{pages}{276} (\bibinfo{year}{2007}).

\bibitem{2000Caspi_PRL}
\bibinfo{author}{Caspi, A.}, \bibinfo{author}{Granek, R.} \&
  \bibinfo{author}{Elbaum, M.}
\newblock \bibinfo{title}{Enhanced diffusion in active intracellular
  transport}.
\newblock \emph{\bibinfo{journal}{Phys. Rev. Lett.}}
  \textbf{\bibinfo{volume}{85}}, \bibinfo{pages}{5655} (\bibinfo{year}{2000}).

\bibitem{2013Becker_PRL}
\bibinfo{author}{Becker, T.}, \bibinfo{author}{Nelissen, K.},
  \bibinfo{author}{Cleuren, B.}, \bibinfo{author}{Partoens, B.} \&
  \bibinfo{author}{Van~den Broeck, C.}
\newblock \bibinfo{title}{Diffusion of interacting particles in discrete
  geometries}.
\newblock \emph{\bibinfo{journal}{Phys. Rev. Lett.}}
  \textbf{\bibinfo{volume}{111}}, \bibinfo{pages}{110601}
  (\bibinfo{year}{2013}).

\bibitem{2016Woodhouse_PNAS}
\bibinfo{author}{Woodhouse, F.~G.}, \bibinfo{author}{Forrow, A.},
  \bibinfo{author}{Fawcett, J.~B.} \& \bibinfo{author}{Dunkel, J.}
\newblock \bibinfo{title}{Stochastic cycle selection in active flow networks}.
\newblock \emph{\bibinfo{journal}{Proc. Natl. Acad. Sci. U.S.A.}}
  \textbf{\bibinfo{volume}{113}}, \bibinfo{pages}{8200--8205}
  (\bibinfo{year}{2016}).

\end{thebibliography}
\end{document}